\documentclass{article}
\usepackage[leqno]{amsmath}
\usepackage{amsthm,amscd}
\usepackage{amsfonts} 
\usepackage{amssymb} 
\newtheorem{theorem}{Theorem}[section]

\newtheorem{lemma}[theorem]{Lemma}
\newtheorem{proposition}[theorem]{Proposition}

\newcommand{\A}{{\mathcal A}}

\newcommand{\bfX}{{\mathbf X}}

\newcommand{\bfP}{{\mathbf P}}

\newcommand{\sD}{{\mathsf D}}

\newcommand{\cH}{{\mathcal H}}

\newcommand{\Der}{{\rm Der}}

\begin{document}
\title{Bases of the contact-order filtration of
derivations of Coxeter arrangements}
\author{{\sc Hiroaki Terao}
\footnote{partially supported by the Grant-in-aid for scientific research
(No. 14340018 and 13874005), the Ministry of Education, Sports,
Science and Technology, Japan}\\
{\small \it Tokyo Metropolitan University, Mathematics Department}\\
{\small \it Minami-Ohsawa, Hachioji, Tokyo 192-0397, Japan}
}
\date{}
\maketitle

\begin{abstract}

\noindent
In a recent paper we constructed a basis
for the contact-order filtration of the module of derivations on 
the orbit space of a finite real reflection
group acting on an $\ell$-dimensional Euclidean space.
%
Recently M. Yoshinaga constructed another basis
for the contact-order filtration in \cite{yos1}. 
In this note we give an explicit formula
relating Yoshinaga's basis to the basis
given in \cite{ter}. The two bases 
turn out to be equal (up to a constant matrix).
\end{abstract}

\bigskip

\setcounter{section}{0}
\setcounter{equation}{0}

\section{The setup and the main result}

Let 
$V$ be an $\ell$-dimensional Euclidean vector space with 
inner product $I$.
Then its dual space $V^{*}$ is 
equipped with the inner product $I^{*} $
which is induced by $I$.
Let $S$  be the symmetric algebra of $V^{*} $ over $\mathbb R$.
Identify $S$ with the algebra of polynomial functions on
$V$.  
Let ${\rm Der}_S$ be the $S$-module of ${\mathbb R}$-linear
derivations of $S$.
When $X_1, \cdots, X_{\ell}$ denote a basis for $V^*$,
the partial derivations 
$\partial_i := \partial / \partial X_{i} $ 
with respect to 
$X_i~~(1\leq i\leq \ell)$
naturally form a basis for ${\rm Der}_S$ over $S$. 
Let $K$  be the field of quotients of $S$
and $\Der_{K} $ be the $K$-vector space of ${\mathbb R}$-linear
derivations of $K$. Then
the partial derivations 
$\partial_i ~~(1\leq i\leq \ell)$
naturally form a basis for ${\rm Der}_K$ over $K$. 

Let $W$ be a
finite irreducible orthogonal reflection group 
(a {\bf Coxeter group}) acting on $V$.
The Coxeter 
group $W$ naturally acts on $V^{*}$, $S$  and
$\Der_{S}$.  
The $W$-invariant subring of $S$
is denoted by  $R$.  
Then it is classically known \cite[V.5.3, Theorem 3]{bou} that there
exist algebraically independent homogeneous polynomials $P_1, \cdots, 
P_{\ell}
\in R$ with $\deg P_{1} \leq\dots\leq \deg P_{\ell}$,
which 
are called {\bf basic invariants},
 such that $R={\mathbb R}[P_1, \cdots, P_{\ell}].$
The {\bf primitive derivation} $D \in {\rm Der}_R$ is characterized
 by
\begin{eqnarray*}
D P_i=\begin{cases}
                                        1 & \text{ for } i=\ell,  \\
                                        0 & \text{ otherwise.}
     \end{cases}
\end{eqnarray*}
Write
$\bfP = (P_{1}, \dots, P_{\ell}  )$.
Let $J(\bfP)$ denote the Jacobian matrix
$[\partial P_{j}/\partial X_{i}  ]_{ij} $.
Define  
$\Delta := \det J(\bfP)$.
Then $\Delta$ is an anti-invariant \cite[V.5.4]{bou}
and
$\Delta^{2} \in R$. 
Let
$$I : \Der_{R} \times \Der_{R} \to \frac{1}{\Delta^2} R$$ 
be the symmetric $R$-bilinear form induced by $I$.
Let
$$\nabla:
{\rm Der}_{R} \times {\rm Der}_{R}  \longrightarrow  
\frac{1}{\Delta^{2}}{\rm Der}_{R}$$
$$({ X}, { Y})  \longmapsto  \bigtriangledown_
{ X}{ Y}$$
be
the {\bf Levi-Civita connection} with respect to $I$. 
Since $\frac{\partial}{\partial P_j} \in \frac{1}{\Delta}\Der_{S}$ 
by Cramer's rule, one can embed $\Der_{R} $ into $\Der_{K}$.
Extend the Levi-Civita connection naturally to 
\[
\nabla : {\rm Der}_{K} \times {\rm Der}_{K}  \longrightarrow  
{\rm Der}_{K}.
\]
Note
$\nabla_{\xi}(\eta) = \sum_{i} \xi(\eta(X_{i} ))\partial_{i}$
because each $I(\partial_{i}, \partial_{j})$ is constant. 

    In \cite{ter}, we introduced derivations
$\xi^{(m)}_1, \dots, \xi^{(m)}_\ell \in \Der_S$
for $m\geq 0$ by
$$
(\xi^{(m)}_1, \dots, \xi^{(m)}_\ell):= 
\begin{cases}
   (\frac{\partial}{\partial P_1}, \dots, \frac{\partial}{\partial P_\ell})
J({\bf P})^T AJ(D^k \left[{\bf X} \right])^{-1} \text{ if } m=2k, \\
   (\frac{\partial}{\partial P_1}, \dots, \frac{\partial}{\partial P_\ell})
J({\bf P})^T AJ(D^k \left[{\bf X} \right])^{-1} J({\bf P}) \text{ if } m=2k+1.
   \end{cases}
$$
Here $A := \left[I^{*} (X_{i}, X_{j})\right]_{ij} $ and
$J(D^{k}[\bfX]) := \left[\partial(D^{k}(X_{j})/\partial X_{i}\right]_{ij}
~~(D^{k} := D\circ D\circ\dots \circ D ~~$($k$ times)).
 Let ${\mathcal A}$ be the Coxeter arrangement
determined by the Coxeter group $W$:
 $\A$ is the set of reflecting hyperplanes.
 Choose for each hyperplane $H \in \cal A$ a linear form $\alpha_H
 \in V^* $ such that $H=\ker (\alpha_H).$
The derivations
$\xi^{(m)}_1, \dots, \xi^{(m)}_\ell \in \Der_S$
were constructed in \cite{ter} so that they may form
a basis for 
the $S$-module \cite{zie}
  \begin{eqnarray*}
 \sD^{(m)}({\mathcal A}):=\{\theta \in {\rm Der}_S ~|~ \theta (\alpha_H) \in 
 S\alpha_{H}^{m} {\rm{~for~ any~}} H \in {\cal A} \}
 \end{eqnarray*}
 for each nonnegative integer $m$.
The filtration
 $$
 {\rm Der}_S = \sD^{(0)}({\mathcal A}) \supset \sD^{(1)}({\mathcal A}) 
 \supset \sD^{(2)}({\mathcal A}) 
 \supset \cdots 
 $$
of $\Der_{S} $ is called
the {\bf contact-order filtration}.
Define
\[
{\mathcal H^{(k)}} =
\sD^{(2k-1)}(\A) \cap \Der_{R}  
\]
for each $k\geq 1$.  Let ${\mathcal H^{(0)}}=\Der_{R}$.
Then we have a filtration
 $$
 {\rm Der}_R = \cH^{(0)} \supset \cH^{(1)}
 \supset \cH^{(2)} 
 \supset \cdots 
 $$
of $\Der_{R}$.  This filtration is known \cite{ter2}
to be equal to the {\bf Hodge filtration} introduced by
K. Saito \cite{sai2}.   
The derivations
$\xi^{(2k-1)}_1, \dots, \xi^{(2k-1)}_\ell \in \Der_R$
form an $R$-basis for $\cH^{(k)} $ 
for $k\geq 1$.  
Therefore the contact-order filtration, retricted to 
$\Der_{R}$, gives the Hodge filtration. 
Define
\[
T := \{f\in R \mid Df = 0\} = {\mathbb R}[P_{1}, \dots ,P_{\ell-1}].
\]
Then the  covariant derivative
$\nabla_{D} : \Der_{K} {\longrightarrow} {\rm Der}_{K}$
is $T$-linear.
The Hodge filtration was originally defined 
so that the $T$-linear map
\[
\nabla_{D}^{k}:=\nabla_{D}\circ\dots\circ\nabla_{D}~(k {\rm ~times} )
:
\cH^{(k)}  \longrightarrow  \cH^{(0)}=\Der_{R}      
\]
is bijective
\cite{sai2}.  Thus we may define
$\nabla_{D}^{-k} \xi \in \cH^{(k)}  $ for any
$\xi\in\Der_{R} $ and $k\geq 0$.
In \cite{yos1} M. Yoshinaga proved the following

\begin{theorem}
\label{yoshinagastheorem} 
{\rm (\cite[Theorem 6]{yos1})}
Let $k\geq 0$ and 
$$E:= \sum_{i} X_{i} \partial_{i} 
=\sum_{i} (\deg P_{i} ) P_{i} (\partial/\partial P_{i})
$$ be the Euler derivation. 
Suppose that
$\xi_{1},\dots, \xi_{\ell}  $ are a basis for $\sD^{(1)}(\A) $.
Then

(1) 
the derivations
$
\nabla_{\xi_{1} } \nabla_{D}^{-k} E,
\dots
,
\nabla_{\xi_{\ell} } \nabla_{D}^{-k} E
$
form a basis for 
$\sD^{(2k-1)}(\A)$ over $S$,
and

(2) the derivations
$
\nabla_{\partial_{1} } \nabla_{D}^{-k} E,
\dots
,
\nabla_{\partial_{\ell} } \nabla_{D}^{-k} E
$
form a basis for 
$\sD^{(2k)} (\A)$ over $S$.

\end{theorem}

It thus seems natural to ask how the basis above constructed 
by Yoshinaga is related to
the basis $\xi^{(m)}_{1}, \dots, \xi^{(m)}_{\ell}$ given in \cite{ter}. 
The following theorem answers this question:

\begin{theorem}
\label{maintheorem} 
Let $k\geq 0$.  Then 
\[
(\xi^{(2k+1)}_{1}, \dots, \xi^{(2k+1)}_{\ell})
=
(-1)^{k} 
(\nabla_{\xi^{(1)} _{1} } \nabla_{D}^{-k} E,
\dots
,
\nabla_{\xi^{(1)} _{\ell} } \nabla_{D}^{-k} E
)
\]
and
\[
(\xi^{(2k)}_{1}, \dots, \xi^{(2k)}_{\ell})
=
(-1)^{k} 
(
\nabla_{\partial_{1} } \nabla_{D}^{-k} E,
\dots
,
\nabla_{\partial_{\ell} } \nabla_{D}^{-k} E
)
A.
\]

\end{theorem}

So the two bases turn out to be equal up to a
constant matrix.

\smallskip

The significance of Theorem \ref{maintheorem} 
is as follows:
Recently M. Yoshinaga \cite{yos2}
affirmatively settled the Edelman-Reiner conjecture
which asserts that the cones of the
extended Shi/Catalan arrangements are 
free.  In other words, there exist basic
derivations for each extended Shi/Catalan arrangement.
However, a formula for 
basic derivations is still unknown.
Since  the derivations 
$\xi^{(m)}_{i}  (i = 1,\dots, \ell)$ 
are the principal (=highest degree) parts of basic
derivations, Theorem \ref{maintheorem} 
can be interpreted as a differential-geometric formula
for the ``principal parts.'' 
So it may suggest the existence of a differential-geometric
formula for the {\it whole} basic
derivations,
including the ``non-principal part.'' 
One may also regard Theorem \ref{maintheorem}
as a very explicit algebraic description of the 
derivations in the right hand side.
They are, as mentioned above, when $m$  is odd, 
bases for the Hodge filtration
which is the key to define the flat structure on
the orbit space $V/W$ in \cite{sai4}.
The flat structure is called the {\it Frobenius 
manifold structure} from the view point of topological field theory
\cite{dub1}.

\section{Proof}

We will prove Theorem \ref{maintheorem} in this section.
First we show the following.

\begin{lemma}
\label{2.1}
For $k\geq 1$  and $\xi\in 
\{
\xi^{(1)}_{1},
\dots
\xi^{(1)}_{\ell}  
\} $,
we have 
\[
\nabla_{D}^{k} \circ\nabla_{\xi} 
-
\nabla_{\xi}\circ\nabla_{D}^{k}  
=
k \,\,\nabla_{D}^{k-1} \circ\nabla_{[D, \xi]}. 
\]
\end{lemma}

\begin{proof} 
We use an induction on $k$.
When $k=1$, the lemma
asserts 
\[
\nabla_{D} \circ\nabla_{\xi} 
-
\nabla_{\xi}\circ\nabla_{D}  
=
\nabla_{[D, \xi]}, 
\]
  which is the integrable property of the 
  Levi-Civita connection $\nabla$.
    Let $k>1$.
  We have
  \begin{eqnarray*}   
\nabla_{D}^{k} \circ\nabla_{\xi} 
&=&
\nabla_{D}^{k-1}\circ(\nabla_{\xi}\circ\nabla_{D}  +
\nabla_{[D, \xi]})\\
&=&
(\nabla_{D}^{k-1}\circ\nabla_{\xi})\circ\nabla_{D}  +
\nabla_{D}^{k-1}\circ\nabla_{[D, \xi]}\\
&=&
(\nabla_{\xi}\circ\nabla_{D}^{k-1}+(k-1)
\nabla_{D}^{k-2}\circ\nabla_{[D, \xi]})\circ\nabla_{D}  +
\nabla_{D}^{k-1}\circ\nabla_{[D, \xi]}\\
&=&
\nabla_{\xi}\circ\nabla_{D}^{k}
+
(k-1)\nabla_{D}^{k-2}\circ\nabla_{[D, \xi]}\circ\nabla_{D}  
+
\nabla_{D}^{k-1}\circ\nabla_{[D, \xi]}
  \end{eqnarray*}   
  by using the induction assumption.
Let
$1\leq i\leq \ell$.
Since
 $\deg \xi(P_{i} ) < 2 (\deg P_{\ell})$,
we have
$$
[D, [D, \xi]](P_{i}) = D^{2} (\xi(P_{i} ))  = 0
$$
and
$[D, [D, \xi]] = 0.$ 
So we obtain
$$
  \nabla_{D}\circ\nabla_{[D, \xi]}
  =
  \nabla_{[D, \xi]}\circ\nabla_{D}.
  $$
  This implies
  \[
\nabla_{D}^{k} \circ\nabla_{\xi} 
=
\nabla_{\xi}\circ\nabla_{D}^{k}
+
k\,\,\nabla_{D}^{k-1}\circ\nabla_{[D, \xi]}.  
  \]
\end{proof}

Recall the $W$-invariant inner product
$I^{*} : V^{*} \times V^{*} \longrightarrow {\mathbb R}$.
Let $\Omega^{1}_{R}$
denote the $R$-module of K\"ahler differentials.
Let 
\[
I^{*} : \Omega^{1}_{R} \times \Omega^{1}_{R} \longrightarrow R
\]
  be the symmetric $R$-bilinear form induced by $I^{*}$.  
Let 
\[ 
G :=[I^*(dP_i, dP_j)]_{ij}= J({\bf P})^T A J({\bf P}).
\]
Define
$$
B^{(k)}:=-J({\bf P})^T AJ(D^k[{\bf X}])J(D^{k-1}[{\bf X}])^{-1}J({\bf P})
$$
for $k \geq 1$ as in \cite{ter}. 
Then we have

\begin{lemma}
\label{2.2} 

(1) Every entry of $B^{(k)}$  lies in $T$;
$D[B^{(k)}]=0$, 

(2) $\det B^{(k)} \in {\mathbb R}^{*}  $, 

(3) $D[G] = B^{(1)} + (B^{(1)})^{T}$, 

(4) $B^{(k+1)}=B^{(1)}+ k D[G].$ 
\end{lemma}

\begin{proof}
By \cite[Lemmas 3.3 - 3.6, Remark 3.7]{ter}. 
\end{proof}

The following proposition is one of the main results in
\cite{ter2}:

\begin{proposition}
\label{2.3}
For $k \geq 1,$ 

$$(\nabla^k_D \xi^{(2k-1)}_1, \dots, \nabla^k_D \xi^{(2k-1)}_\ell)
=(-1)^{k-1}
(\frac{\partial}{\partial P_1}, \dots, \frac{\partial}
{\partial P_\ell}) 
B^{(k)}.$$  
\end{proposition}

In the rest of this note,
let 
$(\xi_1, \dots, \xi_\ell)$ 
denote
$(\xi^{(1)}_1, \dots, \xi^{(1)}_\ell)$ 
for simplicity.
Note
\[ 
(\xi_1, \dots, \xi_\ell)
=
(\frac{\partial}{\partial P_1}, \dots, \frac{\partial}
{\partial P_\ell}) G.
\]

\begin{lemma}
\label{2.4}
$$
([D, \xi_{1}],
\dots
[D, \xi_{\ell}])
=
(\frac{\partial}{\partial P_1}, \dots, \frac{\partial}
{\partial P_\ell}) D[G]
=
(\nabla_{D}\xi_1, \dots, \nabla_{D}\xi_\ell)
(B^{(1)})^{-1} D[G].
$$
\end{lemma}

\begin{proof}
The first equality easily follows from
\[
[D, \xi_{j}](P_{i} )
=
D\circ\xi_{j}(P_{i} )
-
\xi_{j}\circ D(P_{i} )
=
D\circ\xi_{j}(P_{i} )
=
D[I^{*}(dP_{i}, dP_{j}  ) ]. 
\]
For the second equality, apply
Proposition \ref{2.3} when $k=1$.
\end{proof}

\begin{lemma}
\label{2.6}
For any $\eta\in \Der_{K}, $
\[
(
\nabla_{\xi_{1}} \eta,
\dots, 
\nabla_{\xi_{\ell}} \eta)
=
(
\nabla_{\partial_{1}} \eta,
\dots, 
\nabla_{\partial_{\ell}} \eta)
A
J(\bfP).
\]
  \end{lemma}

\begin{proof}
Since both sides are
additive with respect to $\eta$,
we may assume $\eta = f \partial_{j}$ 
for some $j$ and $f\in K$.
Then we have
\begin{eqnarray*}
(
\nabla_{\xi_{1}} \eta,
\dots, 
\nabla_{\xi_{\ell}} \eta)
&=&
(
\nabla_{\xi_{1}} (f \partial_{j}),
\dots, 
\nabla_{\xi_{\ell}} (f \partial_{j})
)\\
&=&
(
\xi_{1} (f)\partial_{j},
\dots, 
\xi_{\ell} (f)\partial_{j})\\
&=&
(
(\partial_{1} f)\partial_{j}
,
\dots
,
(\partial_{\ell} f)\partial_{j}
)
A J(\bfP)\\
&=&
(\nabla_{\partial_{1}} (f \partial_{j}),
\dots, 
\nabla_{\partial_{\ell}} (f \partial_{j}))
A
J(\bfP)\\
&=&
(\nabla_{\partial_{1}} \eta,
\dots, 
\nabla_{\partial_{\ell}} \eta)
A
J(\bfP).
\end{eqnarray*}
\end{proof}

\bigskip
\noindent
{\bf Proof of Theorem \ref{maintheorem}}.
When $k=0$, the formulas are obviously true. 
Let $k\geq 1$.
Let $\xi\in\{\xi_{1}, \dots , \xi_{\ell}\}$.
Apply the identity
in Lemma \ref{2.1} to $\nabla_{D}^{-k} E$ to get
\[
(\nabla_{D}^{k}\circ\nabla_{\xi}) (\nabla_{D}^{-k} E) 
-
\nabla_{\xi} E
=
k 
(
\nabla_{D}^{k-1} 
\circ
\nabla_{[D, \xi]}
)(
\nabla_{D}^{-k} E)
   =
k \,\,
\nabla_{D}^{-1}
(
\nabla_{[D, \xi]} E
).
\]
Since $\nabla_{\eta} E = \eta$ for any $\eta\in\Der_{K},$
one obtains 
\[
\nabla_{D}^{k}
(\nabla_{\xi} \nabla_{D}^{-k} E)
-
  \xi
 =
k \,\,
\nabla_{D}^{-1} [D, \xi]
\]
and thus
\[
\nabla_{D}^{k+1}
(\nabla_{\xi} \nabla_{D}^{-k} E)
=
 \nabla_{D}  \xi
 +
k \,\, [D, \xi]
\]
On the other hand,
Lemma \ref{2.4} asserts
\[
([D, \xi_{1}],
\dots
[D, \xi_{\ell}])
=
(\nabla_{D}\xi_1, \dots, \nabla_{D}
\xi_\ell)
(B^{(1)})^{-1} D[G].
\]
So we have
\begin{eqnarray*}
&~&
(
\nabla_{D}^{k+1}  
(
\nabla_{\xi_{1} }  
\nabla_{D}^{-k} E
)
,
\dots
,
\nabla_{D}^{k+1}  
(
\nabla_{\xi_{\ell} }  \nabla_{D}^{-k} E)
)\\
&=&
(\nabla_{D}\xi_1, \dots, \nabla_{D}\xi_\ell)
+
k\,
(\nabla_{D}\xi_1, \dots, \nabla_{D}\xi_\ell)
(B^{(1)} )^{-1} 
D[G]\\
&=& 
(\nabla_{D}\xi_1, \dots, \nabla_{D}\xi_\ell)
(B^{(1)} )^{-1} 
(B^{(1)} + k D[G])\\
&=& 
(\nabla_{D}\xi_1, \dots, \nabla_{D}\xi_\ell)
(B^{(1)} )^{-1} 
B^{(k+1)}\\
&=& 
(
\partial/\partial P_1
, \dots, 
\partial/\partial P_\ell
)
B^{(k+1)}\\
&=&
(-1)^{k}
(\nabla_{D}^{k+1} \xi^{(2k+1)}_1, \dots, 
\nabla_{D}^{k+1} \xi^{(2k+1)}_\ell)
\end{eqnarray*}
by Lemma \ref{2.2} and Proposition \ref{2.3}.
This proves the first formula.

For the second formula,
compute
\[
(
\nabla_{\xi_{1}} \nabla_{D}^{-k} E,
\dots, 
\nabla_{\xi_{\ell}} \nabla_{D}^{-k} E
)
=
(
\nabla_{\partial_{1}} \nabla_{D}^{-k} E,
\dots, 
\nabla_{\partial_{\ell}} \nabla_{D}^{-k} E
)
A
J(\bfP)
\]
by Lemma \ref{2.6}.
Apply the first formula and we get
\begin{eqnarray*}
(
\xi^{(2k)}_{1},
\dots, 
\xi^{(2k)}_{\ell}
)
J(\bfP)
&=&
(
\xi^{(2k+1)}_{1},
\dots, 
\xi^{(2k+1)}_{\ell}
)
\\
&=&
(-1)^{k} 
(
\nabla_{\xi_{1}} \nabla_{D}^{-k} E,
\dots, 
\nabla_{\xi_{\ell}} \nabla_{D}^{-k} E
)
\\
&=&
(-1)^{k} 
(
\nabla_{\partial_{1}} \nabla_{D}^{-k} E,
\dots, 
\nabla_{\partial_{\ell}} \nabla_{D}^{-k} E
)
A
J(\bfP). 
\end{eqnarray*}
This completes the proof of Theorem \ref{maintheorem}.
\begin{flushright}
$\square$ 
\end{flushright}

\end{document}